\documentclass[11pt]{article}
\usepackage{amsmath,amsfonts,amscd,amssymb,theorem}

 \newcommand{\doubleheaddownarrow}{\big\downarrow\kern-3.325mm\downarrow}

 \newcommand{\PP}{\mathbb P}
 \newcommand{\Q}{\mathbb Q}
 \newcommand{\Z}{\mathbb Z}
 \newcommand{\Oh}{\mathcal O}

 \newcommand{\II}{\operatorname{II}} 
  
 \newcommand{\alt}{\operatorname{alt}}

 \newcommand{\rank}{\operatorname{rank}}
 \newcommand{\ad}{\operatorname{ad}}

 \newcommand{\SL}{\operatorname{SL}}
 
 \newtheorem{theorem}{Theorem}[section]
 \newtheorem{lemma}[theorem]{Lemma}
 \newtheorem{prop}[theorem]{Proposition}
 
 {
 \theorembodyfont{\rmfamily}

 \newtheorem{rem}[theorem]{Remark}

 }
 \newenvironment{pf}{\paragraph{Proof}}{\par\medskip}
 
 \newcommand{\qed}{\ifhmode\unskip\nobreak\fi\quad\ensuremath\square}
\newcommand{\QED}{\ifhmode\unskip\nobreak\fi\quad\ensuremath{\mathrm{QED}}}

\numberwithin{equation}{section}
\title{Some calculations on type $\II_1$ unprojection}
\author{Stavros Argyrios Papadakis}
\date{August 2007}

\begin{document}
\maketitle

\begin {abstract}
The type $\II_1$ unprojection is, by definition, the generic complete intersection type II 
unprojection,  in the sense of \cite{P}~Section~3.1, for the 
parameter value $k = 1$, and depends on a parameter $n \geq 2$. Our main results 
are the explicit calculation of the linear relations  of the type $\II_1$ unprojection 
for any value $n \geq 2$  (Theorem~\ref{thm!aboutlinears}) and the explicit calculation 
of the quadratic equation for the case $n=3$  (Theorem~\ref{thm!quadratic}).
In addition, Section~\ref{sec!applications} contains applications to algebraic geometry,
while Section~\ref{sec!appendix} contains the Macaulay 2 code for the type $\II_1$
unprojection for the parameter value $n=3$.
\end{abstract}

\section {Introduction}   \label{sec!Introduction}

The first appearance of the type II unprojection was in the study of  elliptic 
involutions between Fano $3$-fold hypersurfaces in \cite{CPR}, while the theoretical 
foundations were developed in \cite{P} using valuations. Valuations were useful to prove 
the existence of certain relations, but didn't provide any
explicit formulas for them. The present paper is an effort towards the explicit
calculations for the  type II unprojection.

We define the type $\II_k$ unprojection to be the generic complete 
intersection type II  unprojection,  in the sense of \cite{P}~Section~3.1, for the 
parameter value $k \geq 1$.  According to \cite{P}, it depends on a parameter
$n \geq 2$, increases the codimension from $nk-1$   to  $nk-1+(k+1)$ and preserves
Gorensteiness.
  To our knowledge, the only previously done explicit calculation for the type II
unprojection  was for the type $\II_1$ case for $n=2$ (\cite{CPR}, \cite{R}), which
is reproduced for completeness in Subsection~\ref{sub!reidcasenis2} below.

Section~\ref{sec!linears} contains the calculation, using homological and multilinear
algebra, of the linear equations of the type  $\II_1$ unprojection for any value $n \geq 2$.
The main result is Theorem~\ref{thm!aboutlinears}, which provides explicit formulas
for the linear relations.

Section~\ref{sec!quadratic} contains the calculation of the quadratic equation of the type  
$\II_1$  unprojection for the case $n=3$. The main result is Theorem~\ref{thm!quadratic}
which provides an explicit symmetric -- in the sense of $\SL_3$ invariance, see
Subsection~\ref{subs!invariance} -- formula for the quadratic relation for this case.
The proof of Theorem~\ref{thm!quadratic} is based on the explicit equality 
(\ref{eqn!bigcheck}) which was verified using the computer algebra program Macaulay 2
\cite{GS93-08}. In Subsection~\ref{sub!howwegot} we briefly sketch how we arrived to the 
formulas contained in Theorem~\ref{thm!quadratic} by computer assisted calculations (using
the computer algebra program Maple).

As an application of the above results, we sketch in Section~\ref{sec!applications}  
the construction of two codimension $4$ Fano $3$-folds. The quasismoothness checking,
using the explicit equations obtained in Sections~\ref{sec!linears} and \ref{sec!quadratic},
was done by the computer algebra program Singular~\cite{GPS01}.

The calculation of the quadratic equation for the type $\II_1$ unprojection
for $n \geq 4$ remains open. We believe that the explicit formulas of
the linear relations obtained in the present work together with representation
and invariant theoretic techniques could, perhaps, lead to their calculation,
and also to a better understanding of the complicated looking formula 
(\ref{eqn!dfnofb}).  To our knowledge, the problem
of calculating the type $\II_k$ unprojection for $k \geq 2$ is also open.

ACKNOWLEDGEMENTS: I wish to thank Janko Boehm, Stephen Donkin, S\'ebastien Jansou,
Nondas Kechagias, Frank--Olaf Schreyer and Bart Van Steirteghem  for  useful discussions.
Parts of this work were financially supported by 
the Deutschen Forschungsgemeinschaft Schr 307/4-2 and by the Portuguese 
Funda\c{c}\~ao para a Ci\^encia e a Tecno\-lo\-gia through Grant \\
SFRH/BPD/22846/2005.

\section {Notation}   \label{sec!notation}

As already mentioned in Section~\ref{sec!Introduction} by 
type $\II_1$ unprojection we mean the
generic complete intersection type II unprojection
in the sense of \cite{P}~Section~3.1 with fixed parameter value $k=1$.
According to \cite{P}, the type $\II_1$ unprojection depends 
on an integral parameter $n \geq 2$, the initial data for the unprojection
is specified by a triple 
\begin {equation}  \label {eq!thetriple}
   I_X  \subset I_D  \subset \Oh_{amb}
\end {equation}
(where $I_X$ and $I_D$ are ideals of  $\Oh_{amb}$), and
the unprojection constructs an ideal 
\begin {equation}  \label {eq!thefinalproduct}
          I_Y \subset \Oh_{amb}[s_0,s_1],
\end {equation}
where $s_0,s_1$ are new variables.

Fix $n \geq 2$. The ambient ring  $\Oh_{amb}$ will be the polynomial ring 
\[
      \Oh_{amb} = \Z [x_1, \dots , x_n, y_1, \dots ,y_n,z, A^p_{ij},B^p_{lm} ] 
\]
where $1 \leq p \leq n-1$,  $1 \leq i < j \leq n$ and
$1 \leq l \leq m \leq n$. The ring $\Oh_{amb}$ corresponds to the ring 
denoted by  $\Oh_{amb}^g$  in \cite {P}~Section~3.1, the variable 
$x_i$ corresponds to the  variable $a_{1i}$ in \cite {P}~Section~2, and
the variable $y_i$   corresponds to  the variable $a_{2i}$ in \cite {P}~Section~2.
(In Section~\ref{sec!quadratic}  $\; \Oh_{amb}$ will change slightly to allow 
$2$ to be invertible.)

We denote by $M$ the $2 \times 2n$ matrix
\begin{equation}     \label{eq!dfnofM}
    M =  \begin {pmatrix}
             y_1 & \dots & y_n  &  z x_1 & \dots & z x_n  \\
             x_1 & \dots & x_n  &  y_1   & \dots & y_n 
     \end {pmatrix} 
\end{equation} 
which corresponds  to the matrix with the  same name in 
\cite{P} Equation~(2.1).

The ideal $I_D \subset \Oh_{amb}$ in (\ref{eq!thetriple})  
is the ideal generated
by the $2 \times 2$ minors of $M$ and corresponds to the
ideal $I_D^g$  in \cite {P}~Section~3.1.

For $1 \leq  p \leq n-1$, we denote by  $A^p$  the $n\times n$ 
skew--symmetric matrix with $(ij)$-entry (for $i <j$)  equal to
$A_{ij}^p$ (and zero diagonal entries),  and by $B^p$ the 
$n\times n$  symmetric 
matrix with $(ij)$-entry (for $i \leq j$) equal to $B^p_{ij}$.

For $1 \leq p \leq n-1$ we set
\begin{equation*}  
   f_p = \sum_{1\leq i < j \leq n }A_{ij}^pF_{ij} + 
   \sum_{i=1}^n B_{ii}^pG_{ii}  +  
      \sum_{1\leq i < j \leq n }B_{ij}^pG_{ij}, 
\end{equation*}
where, 
\[
  F_{ij} = x_i y_j - x_j y_i, \quad  \quad G_{ij} = 2(y_iy_j - z x_ix_j), \quad \quad G_{mm} = y_m^2 - z x_m^2
\]
for $1 \leq i < j \leq n \;$ and  $\;1 \leq m \leq n$. In matrix notation
\begin{equation}  \label{eq!dfnoff}
   f_p = xA^py^t + yB^py^t - zxB^px^t.
\end{equation}
The ideal $I_X$ in (\ref{eq!thetriple})  is the ideal
\[
   I_X = (f_1, \dots ,f_{n-1})
\]
and corresponds to the ideal $I_X^g$ in \cite {P}~Section~3.1.

In addition, we set  
\begin{equation}   \label{eqn!defnofxandy}
     x=(x_1,\dots ,x_n),  \quad y=(y_1, \dots ,y_n).
\end{equation}

\begin {rem} \label{rem!aboutIY}
According to \cite{P}~Proposition~2.16,
the ideal $I_Y$ of  $\Oh_{amb}[s_0,s_1]$ contains $I_X$ as a subset and 
is generated by $I_X$ together with $n$  polynomials
$l_1, \dots ,l_n \in \Oh_{amb}[s_0,s_1]$ 
affine linear in $s_0,s_1$ of the form
\[
   l_i     =  zx_i s_0 + y_i s_1 +\sigma_i^1,
\]
where $1 \leq i \leq n$ and $\sigma_i^1 \in \Oh_{amb}$, 
together with $n$  polynomials  
$l_{n+1}, \dots ,l_{2n} \in \Oh_{amb}[s_0,s_1]$ 
affine linear in $s_0,s_1$ of the form
\[
       l_{i+n} =  y_i s_0 + x_i s_1 + \sigma_i^2 ,
\]
where $1 \leq i \leq n$ and $\sigma_i^2 \in \Oh_{amb}$, 
together with a single affine quadratic polynomial $q \in \Oh_{amb}[s_0,s_1]$
of the form 
\[
     q = s_1^2 - zs_0^2 + \text { lower terms, }
\]
where 'lower terms' means an affine linear polynomial in
$s_0,s_1$, i.e., of the form $e_1s_0+e_2s_1+e_3$ with 
$e_1,e_2,e_3 \in \Oh_{amb}$.

In other words, we have the equality 
\begin {equation} \label {eqn!formatofIY}
   I_Y = (f_1, \dots ,f_{n-1}) + (l_1, \dots , l_{2n}) + (q)
\end {equation}
of ideals of $\Oh_{amb}[s_0,s_1]$. 
Calculating $I_Y$ means providing explicit 
formulas for $l_1,l_2, \dots , l_{2n}$ and $q$.
\end {rem}

For simplicity,  we call the elements of $I_X$ the original relations,
we call $l_1, \dots ,l_{2n}$ the linear relations and, finally, we call 
$q$ the quadratic relation.

For the parameter value $n=2$ the ideal $I_Y$ has been explicitly
calculated by Reid in \cite{R}~Section~9.5. For completeness, we
reproduce his results in Subsection~\ref{sub!reidcasenis2}.

In Theorem~\ref{thm!aboutlinears} we calculate, for any $n \geq 2$, the linear
relations  $l_1, \dots ,l_{2n}$. In addition, in Theorem~\ref{thm!quadratic} we 
calculate the quadratic polynomial $q$ for the parameter value $n=3$.

The calculation of the quadratic polynomial $q$ for $n \geq 4$ remains open.

\section {Linear relations of Type $\II_1$ unprojection} \label{sec!linears}

In this section we use the notations of  Section~\ref{sec!notation}.

\subsection {Generalities}  \label{subs!Koszulcmplx}

Assume $R$ is a ring, $L$ is an $R$-module, and $h \colon L \to R$ a
homomorphism of $R$-modules. For $p \geq 1$, \cite{BH}~p.~43 defines an
$R$-homomorphism $d^ph : \wedge^p L \to \wedge^{p-1} L$  by
\begin {equation}    \label{eq!dfnofkoszul}
   d^ph (v_1 \wedge  \dots \wedge v_p) = \sum_{i=1}^p (-1)^{i+1}
          h(v_i) v_1 \wedge \dots \wedge   \widehat{v}_i 
   \wedge \dots \wedge v_p,
\end{equation}
for all  $v_1, \dots ,v_p \in L$. (For $p=1$ we set $\wedge^0L = R$
and $d^1h = h$.)
The maps $d^ph$  define the Koszul complex 
\[
    \dots \to \wedge^pL \xrightarrow{d^ph}  \wedge^{p-1}L \to \dots 
       \to  \wedge^2L  \xrightarrow{d^2h} L 
        \xrightarrow{h}  R \to0
\]
associated to the homomorphism $h$.

\begin {prop} \label {prop!commutativity}
Assume $h_1,h_2 \colon L \to R$ are two homomorphisms of $R$-modules. 
For $p \geq 2$ we have the equality
\[
    d^{p-1}h_1 \circ d^ph_2 + d^{p-1}h_2 \circ d^ph_1 = 0,
\]
of maps $\wedge^pL \to \wedge^{p-2}L$.
\end {prop}

\begin {pf} 
Indeed,
\[   0 = d^{p-1} (h_1+h_2) \circ d^{p} (h_1+h_2) =   
     d^{p-1}h_1 \circ d^ph_2 + d^{p-1}h_2 \circ d^pp_1, 
\]
since
\[   
   d^{p-1} h_1 \circ d^{p}h_1 =  d^{p-1} h_2 \circ d^{p} h_2 = 0.
\]
\medskip \QED \end{pf}

\subsection {The second complex}    \label{subs!secondcmplx}

Assume $L$ is a free $\Oh_{amb}$-module of rank $n$, and let $e_1, \dots ,e_n$
be  a fixed basis of $L$. We define four $\Oh_{amb}$-homomorphisms
$h_1, \dots ,h_4 : L \to \Oh_{amb}$, with
\[
   h_1(e_i) = y_i, \ \ h_2(e_i) = x_i, \ \ h_3(e_i) = zx_i, \ \ 
  h_4(e_i) = y_i,
\]
for $1 \leq i \leq n$. In addition, we define, for $1 \leq p \leq n$, homomorphisms
\[ \phi_p \colon 
  \wedge^p L \oplus \wedge^p L \to \wedge^{p-1} L \oplus \wedge^{p-1}  L
\]
(where by definition  $\wedge^{0} L = \Oh_{amb}$)
as follows: 
\[ 
  \phi_p (a,b) =
       \begin {pmatrix}
             d^ph_1 & d^ph_3   \\ 
             d^ph_2 & d^ph_4 
   \end{pmatrix}  
       \begin {pmatrix}
             a   \\    
             b
   \end{pmatrix} = (d^ph_1(a) + d^ph_3(b), d^ph_2 (a)+ d^ph_4(b)),
\]
for $a,b \in \wedge^p L$.

\begin {prop}
For $p \geq 2$ we have 
 \[  
     \phi_{p-1} \circ  \phi_{p} = 0. 
\]
\end{prop}

\begin {pf}
  Using the equalities  $h_4 = h_1$ and $h_3 = zh_2$
   the result follows  from Proposition~\ref{prop!commutativity}.
\medskip \QED \end {pf}  

The proof of the following proposition will be given in the Subsection~\ref{sub!proofofexactness}.

\begin {prop}  \label{prop!proofofexactness}
The complex 
\begin {equation}   \label{prop!dfnofL}
    {\bf L^*:} \ \ \   \Oh_{amb} \oplus \Oh_{amb} \xleftarrow{\phi_1} L \oplus L  
       \xleftarrow{\phi_2}  \wedge^2 L \oplus \wedge^2 L
      \xleftarrow{\phi_3}  \wedge^3 L  \oplus \wedge^3 L \leftarrow \dots
\end{equation} 
is exact.
\end{prop}

\subsection {Proof of Proposition~\ref{prop!proofofexactness}}  \label {sub!proofofexactness}

The proof of Proposition~\ref{prop!proofofexactness} will be based on the Buchsbaum--Eisenbud
acyclicity criterion as stated in \cite{BH}~Theorem~1.4.13. We first need the following 
combinatorial lemma.

\begin {lemma}  \label{lem!binomialidentity}
Let $p$ be an integer with  $1 \leq p \leq n$. Then 
\[
   \sum_{j=p}^{n} (-1)^{j-p}\binom{n}{j} = \binom {n-1}{p-1}.
\]
\end {lemma}

\begin {pf}  
   We have 
\begin {eqnarray*}
   \sum_{j=p}^{n} (-1)^{j-p}\binom{n}{j} = (-1)^{n-p} +  \sum_{j=p}^{n-1} (-1)^{j-p}\binom{n}{j} \phantom{=======} \\
     \phantom {=====}  = (-1)^{n-p} +  \sum_{j=p}^{n-1} (-1)^{j-p}( \binom{n-1}{j-1}+\binom{n-1}{j})  \\
      \phantom {=====}   = (-1)^{n-p} + \binom{n-1}{p-1} + (-1)^{n-1-p}  =  \binom{n-1}{p-1},
\end {eqnarray*}
where we used the Pascal's rule 
  $ \binom{n}{j}  =  \binom{n-1}{j-1}+\binom{n-1}{j}$  which is valid
whenever $1 \leq j \leq n-1$. \medskip  \QED \end {pf}

For the following we fix an integer $p$ with $1 \leq p \leq n$.  We set 
\[
       r_p =  \sum_{j=p}^{n} (-1)^{j-p} \rank \;(\wedge^p L \oplus \wedge^p L),
\]
where for a finitely generated free $\Oh_{amb}$-module $N$
we denote by $\rank N$ the minimal number of generators of $N$ as $\Oh_{amb}$-module.
Using Lemma~\ref{lem!binomialidentity} we get  
\begin {equation}  \label {eqn!forrp}
      r_p = 2\binom {n-1}{p-1}.
\end {equation}

Denote by  $I_{r_p}(\phi_p)$ the $r_p$-th
Fitting ideal of the map $\phi_p$, that is the
ideal of $\Oh_{amb}$ generated by the $r_p \times r_p$ minors
of any matrix representation of $\phi_p$. (For more information about 
Fitting ideals see, for example, \cite{BH}~Section~1.4.)

We will need the following general property of Koszul complexes.

\begin {lemma}  \label{lemma!diagonalityforkoszul}
Fix $t \in \{1,2\}$  and $s,p$ with  $1 \leq s,p \leq n$.
Denote by $M^p_t$ the matrix representation 
of the map 
\[
        d^{p}h_t \colon \wedge^{p} L  \to \wedge^{p-1} L
\]
with respect to the bases $\; e_{i_1} \wedge \dots \wedge  e_{i_{p}}$ of  $\wedge^p L$ 
and $\;e_{i_1} \wedge \dots \wedge  e_{i_{p-1}}$ of  $\wedge^{p-1} L$ induced by
the basis $(e_1, \dots , e_n)$ of $L$.
For a suitable ordering of the two bases  there exists an 
$(r_p/2) \times (r_p/2)$ submatrix of $M^p_t$ which is diagonal with
diagonal entries equal to  $a_t$ or $-a_t$, where  $a_1 =y_s$
and $a_2=x_s$.

\end {lemma}

\begin {pf}  
Consider the subset of the above mentioned basis  of  $\wedge^p L$ 
consisting of the elements where $e_s$ appears on the wedge. This subset 
has  $\binom{n-1}{p-1}$ elements, hence $(r_p/2)$ elements using 
(\ref{eqn!forrp}). The result follows by 
restricting the homomorphism $d^ph_t$ to
the $\Oh_{amb}$-submodule of  $\wedge^p L$ generated by this subset and
noticing that the corresponding matrix has a diagonal submatrix of
the desired form. \medskip   \QED \end {pf}

\begin {prop}  \label{prop!specialelements}
Fix $s,p$ with with $1 \leq s,p \leq n$.  There exists a polynomial
$g_s \in I_{r_p}(\phi_p)$ of the form  
\begin{equation}  \label{eqn!polynsg1i}
    g_s = y_s^{r_p} + g^1_s,
\end{equation} 
where $g^1_s$ involves only the variables $y_s,z,x_s$
and has degree in $y_s$ strictly less than $r_p$.
\end {prop} 

\begin {pf}
It follows immediately from the definition of $\phi_p$
and Lemma~\ref{lemma!diagonalityforkoszul}. \medskip \QED \end {pf}

\begin {prop}  \label{prop!overafield}
Let $F$ be either the field $\Q$ of rational numbers or a finite
field $\Z/(q)$, where $q$ is a prime. 
The ideal  $\widetilde{I}$  of  $\Oh_{amb} \otimes_\Z F$ generated by  
the image of $I_{r_p}(\phi_p)$  under the natural map 
$\Oh_{amb} \to \Oh_{amb} \otimes_\Z F$ has codimension greater or equal than $n$.
\end {prop}

\begin {pf}
Using Proposition~\ref{prop!specialelements}
it is clear that there exists a monomial ordering for the variables
$y_j, z, x_j$ (with indices $1 \leq j \leq n$) such that the leading term
of  the polynomial $g_s$ appearing in (\ref{eqn!polynsg1i}) is equal 
to $y_s^{r_p}$, for $1 \leq s \leq n$. Since over the field $F$  an ideal  $\widetilde{I}$ 
and the leading term ideal of $\widetilde{I}$  have the  same codimension the 
result follows. \medskip  \QED \end {pf}

Recall (cf.~\cite{BH}~Section~1.2) that the grade of a proper ideal 
$I \subset \Oh_{amb}$  is the length of a maximal $\Oh_{amb}$-regular
sequence contained in $I$.

\begin {prop}  \label{prop!overtheintegers}
The ideal $I_{r_p}(\phi_p)$ has grade  greater or equal than $n$.
\end {prop}

\begin {pf}
Combining Proposition~\ref{prop!overafield} and \cite{L}~Theorem~3.12
we get that $I_{r_p}(\phi_p)$ has codimension in $\Oh_{amb}$  greater or equal than $n$.
Being a polynomial ring over the integers, the ring $\Oh_{amb}$ is Cohen--Macaulay.
Using \cite{BH}~Corollary~2.14 the result follows. \medskip   \QED \end {pf}

We now finish the proof of Proposition~\ref{prop!proofofexactness}.  
Combining Lemma~\ref{lem!binomialidentity} and Proposition~\ref{prop!overtheintegers},
Proposition~\ref{prop!proofofexactness} follows from the
Buchsbaum--Eisenbud acyclicity criterion as stated in 
\cite{BH}~Theorem~1.4.13.

\subsection {The first complex}   \label{subs!firstcmplx}

Fix $n \geq 2 $. Let $N$ be a free $\Oh_{amb}$-module of rank $n-1$, with basis 
$\widetilde{e}_1, \dots , \widetilde{e}_{n-1}$.

We define an $\Oh_{amb}$-homomorphism $\psi : N \to \Oh_{amb}$, with 
$\psi (\widetilde{e}_p)  =f_p $,  for $1 \leq p \leq n-1$, where $f_p$ 
was defined in  (\ref{eq!dfnoff}).
As in Subsection~\ref{subs!Koszulcmplx}, we have the Koszul complex
\begin{equation}
  {\bf N^*:} \ \ \ \   \  \Oh_{amb} \xleftarrow{\psi} N
       \xleftarrow{d^2\psi}  \wedge^2 N 
      \xleftarrow{d^3\psi}  \wedge^3 N   \leftarrow \dots
\end{equation}

\begin {prop}
The sequence $f_1, \dots ,f_{n-1}$ is an $\Oh_{amb}$-regular sequence. As a consequence,
the complex  ${\bf N^*}$ is exact.
\end{prop}
  
\begin{pf}
For an element $g$ of $\Oh_{amb}$ we denote by $\eta (g)$ the
result of substituting to $g$ zero for $z$, zero for 
$A^p_{ij}$ with $1 \leq p \leq n-1$ and $1 \leq i < j \leq n$,
and zero for $B^p_{ij}$ for $1 \leq p \leq n-1$ and 
$1 \leq i \leq j \leq n$  with $(i,j) \not= (p,p)$.
That is, we set zero the variable $z$ and all possible $A^p_{ij}$
and all possible $B^p_{ij}$ with the exception of $B^p_{pp}$ 
(for $1 \leq p \leq n-1$).

It is clear  that for  $1 \leq p \leq n-1$ we have
\[
       \eta(f_p) = B^p_{pp}(y_p)^2.
\]
Therefore, the ideal $ (\eta(f_1), \dots ,\eta(f_{n-1})) \subset \Oh_{amb} $ 
has codimension in $\Oh_{amb}$ equal to $n-1$. As a consequence, the ideal 
$(f_1, \dots ,f_{n-1})$ of  $\Oh_{amb}$ has also codimension $n-1$. Since
$\Oh_{amb}$ is Cohen--Macaulay the sequence $f_1, \dots ,f_{n-1}$ is an 
$\Oh_{amb}$-regular sequence.
The exactness of the Koszul complex  ${\bf N^*}$ follows from 
\cite{BH}~Corollary~1.6.14. \medskip \QED  \end{pf}

\subsection {The first commutative square} \label{sub!firstsquare}

We define two $\Oh_{amb}$-homomorphisms $u_1,u_2 \colon N \to L$, that  will
make the following diagram
\begin {equation*}
\begin{CD}
               \Oh_{amb}  @< \psi << N  \\
                  @V \begin {pmatrix} 1 \\ 0 \end{pmatrix}
   VV   @V   \begin {pmatrix} u_1 \\ u_2 \end{pmatrix} VV   \\
               \Oh_{amb} \oplus \Oh_{amb} @< \phi <<  L \oplus L  
\end{CD}
\end{equation*}
commutative. That is,
 they will satisfy the two conditions
\begin {eqnarray}  
  \psi & = & h_1 \circ u_1 + z h_2 \circ u_2  \label{eq!commutativity1} \\
  0  & = &  h_1 \circ u_2 +  h_2 \circ u_1    \label{eq!commutativity2}
\end{eqnarray}

We set,  for $1 \leq j \leq m$,
\[
   u_1 (\widetilde{e}_j) = \sum_{i=1}^{n}c_{ij}e_i \quad  \quad \quad \text{ and } \quad \quad \quad
   u_2 (\widetilde{e}_j) = \sum_{i=1}^{n}c'_{ij}e_i,
\]
with 
\[
  c_{ij} =  - \sum_{l=1}^{i-1}  A^j_{li}x_l
            + \sum_{l=i+1}^n    A^j_{il}x_l  
            + \sum_{l=1}^i  B^j_{li}y_l
            + \sum_{l=i+1}^n  B^j_{il}y_l,
\]
and 
\[
  c'_{ij} =    - \sum_{l=1}^i   B^j_{li}x_l
               - \sum_{l=i+1}^n B^j_{il}x_l
\]
for $1 \leq i \leq n$ and $1 \leq j \leq n-1$.

In matrix notation, $u_1$ and $u_2$ are given, with respect to
the above defined bases $(\widetilde{e}_j)$ of $N$ and $(e_i)$ of $L$,
by the following two $n \times (n-1)$ matrices:
\begin{eqnarray}
   u_1 &  = & (-A^1x^t+B^1y^t \;|  -A^2x^t+B^2y^t  \;| \dots | -A^{n-1}x^t+B^{n-1}y^t)  \notag  \\
   u_2 & = & (-B^1x^t  \;|  -B^2x^t  \;| \dots | -B^{n-1}x^t).   \label{eqn!matrixformat}
\end{eqnarray}
That is, for $1 \leq j \leq n-1$, 
 $u_1$ is the matrix with $j$-th column equal to   
$-A^jx^t+B^jy^t$, while  $u_2$ is the matrix with $j$-th column equal to
$-B^jx^t$.

\begin {prop}
The maps $u_1$ and $u_2$ satisfy the equations 
(\ref{eq!commutativity1}) and (\ref{eq!commutativity2}).
\end{prop}

\begin {pf}
 An easy calculation using (\ref{eqn!matrixformat}).
\medskip \QED \end{pf}

\subsection {The basic formula}  \label{sub!basicformula}

The key result of the present subsection is 
Proposition~\ref{prop!basicformula},
essentially a formal consequence of (\ref{eq!commutativity2}), 
which relates 
the differentials and the alternating products of the homomorphisms
$h_i$ and $u_j$.
We need the following definitions. 

Assume $a_1,a_2 \colon N \to L$ are two $\Oh_{amb}$-homomorphisms.

For $p \geq 0$ we have a natural induced homomorphism 
$\wedge^p a_1 \colon \wedge^p N \to \wedge^p L$.
By definition, $\wedge^0 a_1$  is the identity map $\Oh_{amb} \to \Oh_{amb}$, 
For $p \geq 1$, the map $\wedge^p a_1$ is uniquely specified by  the property
\[
    \wedge^p a_1 (c_1 \wedge \dots \wedge c_p) = a_1(c_1) \wedge \dots \wedge a_1(c_p)
\]
for all $c_1, \dots ,c_p \in N$. We denote the map $\wedge^p a_1$  also by
$ \alt (a_1^p)$.

We will now define for $p,q \geq 1$ an $\Oh_{amb}$-homomorphism 
\[ 
  \alt (a_1^p, a_2^q) \colon \wedge^{p+q} N  \to  \wedge^{p+q} L.
\]
We set
\[
  \alt (a_1^p,a_2^q) = \sum w_I, 
\]
where the sum is over all subsets $I \subseteq \{1,\dots ,p+q \}$,
with $|I| = p$ and, by definition, $ w_I = w_1 \wedge \dots \wedge w_{p+q}$ with 
$w_i = a_1$ if $i \in I$ and $w_i = a_2$ if $i \notin I$.
More precisely, even though the tensor product map 
\[ 
  w_1 \otimes \dots \otimes w_{p+q} \colon N^{\otimes(p+q)} \to  L^{\otimes(p+q)}
\]
does not induce a map   $\wedge^{p+q} N  \to  \wedge^{p+q} L$,
taking the sum  of those maps  over the family of subsets $I$ as above 
induces a well--defined  map 
$\wedge^{p+q} N  \to  \wedge^{p+q} L$, and this is the map 
denoted by  $\alt (a_1^p,a_2^q)$.

\begin {rem} For $p=q=1$ we have
\[
   \alt(a_1,a_2) (c_1 \wedge c_2) = a_1(c_1) \wedge a_2 (c_2) + a_2(c_1) \wedge a_2(c_2),
\]
for all $c_1,c_2 \in N$. We will also denote  $\alt(a_1,a_2)$ by $a_1 \wedge a_2 + a_2 \wedge a_1$. \\
For  $p=2, q=1$ we have 
\[
  \alt(a_1^2,a_2) = a_1 \wedge a_1 \wedge a_2 +  a_1 \wedge a_2 \wedge a_1 +  a_2 \wedge a_1 \wedge a_1,
\]
in the sense that
\begin {eqnarray*} 
  & & \alt(a_1^2,a_2) (c_1 \wedge c_2 \wedge c_3) =
      a_1(c_1)\wedge a_1(c_2) \wedge a_2(c_3) +  \\
   &  &  \phantom {====} a_1(c_1) \wedge a_2(c_2) \wedge a_1(c_3) +  
             a_2(c_1) \wedge a_1(c_2) \wedge a_1(c_3)
\end {eqnarray*}
for all $c_1,c_2, c_3 \in N$. 
\end {rem}

For the following proposition, which will be proved in  Subsection~\ref{subs!pfofbasicformula},
recall that $u_1,u_2  \colon N \to L$ were 
defined in  Subsection~\ref{sub!firstsquare}, and $h_1,h_2 \colon L \to \Oh_{amb}$
were defined in Subsection~\ref{subs!secondcmplx}. We 
also use the notational conventions $\alt (u_1^{p},u_2^0) =  \alt (u_1^p)$,
$\alt (u_1^0,u_2^q) =  \alt (u_2^q)$,  $\alt(u_1^{p-2},u_2^q) = 0$  for $p=1$ and $q \geq 1$,
and  $\alt(u_1^{p},u_2^{q-2}) = 0$  for $q= 1$ and $p \geq 1$.

\begin {prop}  \label{prop!basicformula}
a) Let $u \colon M \to L$ and  $h \colon L \to \Oh_{amb}$ be two (arbitrary) homomorphisms
of $\Oh_{amb}$-modules. Assume $p \geq 2$. We have
\[
   d^{p-1} h \circ \alt (u^{p-1}) = \alt (u^{p-2}) \circ d^{p-1}(h \circ u).
\] 

b) For $p,q \geq 1$ we have 
\begin{eqnarray}   \label{eq!basicformula}
   (d^{p+q-1}h_2)\circ \alt(u_1^p,u_2^{q-1}) + (d^{p+q-1}h_1) \circ
\alt(u_1^{p-1},u_2^{q})  =  \phantom{==========}   \\
    \alt(u_1^{p-2},u_2^q) \circ
   d^{p+q-1}(h_1 \circ u_1)  + \alt(u_1^{p},u_2^{q-2}) \circ
   d^{p+q-1}(h_2 \circ u_2),  \notag   
\end{eqnarray}
where the equality is, of course, as maps $\quad \wedge^{p+q-1}N \to \wedge^{p+q-2}L$.
\end{prop}

\subsection {Proof of Proposition~\ref{prop!basicformula}} \label{subs!pfofbasicformula}

We will use the following two lemmas.

\begin {lemma}   \label {lem!firstgoodidentity}
Assume $p,q \geq 1$ and $u_1, u_2 \colon N \to L$ 
are two $\Oh_{amb}$-module homomorphisms. We have
\[
   \alt(u_1^p, u_2^q) = u_1 \wedge  \alt(u_1^{p-1}, u_2^q) +
          u_2 \wedge  \alt(u_1^p, u_2^{q-1}),
\]
in the sense that 
\begin {eqnarray*}
   \alt(u_1^p, u_2^q) (c \wedge c_1) = \qquad \qquad  \qquad  \qquad  \qquad  \qquad  \qquad   \qquad  \qquad  \qquad \\
      u_1(c) \wedge   \Bigl[ \alt(u_1^{p-1}, u_2^q)(c_1) \Bigr]  + 
      u_2 (c) \wedge  \Bigl[ \alt(u_1^p, u_2^{q-1})(c_1) \Bigr] 
\end {eqnarray*} 
for all $c \in N$ and $c_1 \in \wedge^{p+q-1}N$.     
\end {lemma}

\begin {pf}  Immediate from the definitions.  \medskip \QED 
\end {pf}

\begin {lemma}   \label {lem!secondgoodidentity}
Assume $p,q \geq 1$, $u_1, u_2 \colon N \to L$ 
are two $\Oh_{amb}$-module homomorphisms, and
$h \colon N \to \Oh_{amb}$ is an $\Oh_{amb}$-homomorphism. We have
\begin {eqnarray*}
   \alt(u_1^p, u_2^q) \circ (d^{p+q+1}h) =    h \wedge \alt(u_1^p,u_2^q) -    
                         u_1 \wedge   \Bigl[\alt(u_1^{p-1}, u_2^q) \circ d^{p+q}h \Bigr] \qquad  \\
                        - u_2 \wedge  \Bigl[\alt(u_1^{p}, u_2^{q-1}) \circ d^{p+q}h \Bigr], \qquad 
\end {eqnarray*} 
in the sense that 
\begin {eqnarray*}
 \Bigl[ \alt(u_1^p, u_2^q) \circ (d^{p+q+1}h)  \Bigr] (c \wedge c_1) =   \qquad  \qquad  \qquad  
    \qquad  \qquad  \qquad  \qquad   \qquad   \\
    h(c) \Bigl[ \alt(u_1^p,u_2^q)(c_1) \Bigr]  - u_1(c) \wedge  (\alt(u_1^{p-1}, u_2^q) \circ d^{p+q}h)(c_1) \qquad  \\
          -  u_2(c) \wedge  (\alt(u_1^{p}, u_2^{q-1}) \circ d^{p+q}h)(c_1) \qquad  
\end {eqnarray*} 
for all $c \in N$ and $c_1 \in \wedge^{p+q}N$.   
\end {lemma}

\begin {pf} 
   Since 
\begin{equation}  \label {eqn!onemore}
   d^{p+q+1}h  (c \wedge c_1)  =  h(c)c_1 - c \wedge  d^{p+q}h (c_1),
\end{equation}
we have, using Lemma~\ref{lem!firstgoodidentity}, that
\begin {eqnarray*}
    \alt(u_1^p, u_2^q) \circ (d^{p+q+1}h)  (c \wedge c_1) =   h(c) (\alt(u_1^p,u_2^q)(c_1))   
    \qquad  \qquad  \qquad     \\
  \qquad  \qquad      \qquad  \qquad - \Bigl[ u_1 \wedge \alt(u_1^{p-1}, u_2^q)+ u_2 \wedge \alt(u_1^p, u_2^{q-1}) \Bigr]
            (c \wedge  d^{p+q}h (c_1))
\end {eqnarray*}
and the result follows.       \medskip \QED 
\end {pf}

We will now prove part a) of Proposition~\ref{prop!basicformula}. Assume 
$p \geq 2$. For $c_1, \dots ,c_{p-1} \in N$ we have 
\begin {eqnarray*}
   d^{p-1} h \circ \alt (u^{p-1}) (c_1 \wedge \dots \wedge c_{p-1}) = 
     d^{p-1} h  (u(c_1) \wedge \dots \wedge u(c_{p-1}))  \\
     =\sum_{i=1}^{p-1} (-1)^{i-1} h(u(c_i)) 
       (u(c_1) \wedge \dots \wedge \widehat{u(c_i)}  \dots  \wedge u(c_{p-1}))  
    \qquad   \\
   =  \alt (u^{p-2}) \circ d^{p-1}(h \circ u)  (c_1 \wedge \dots \wedge c_{p-1}).  
        \qquad  \qquad  \qquad 
\end {eqnarray*}

We will now prove  part b)  of Proposition~\ref{prop!basicformula}.
Assume that $p,q \geq 1$, that $N,L$ are $\Oh_{amb}$-modules
and that $u_1, u_2 \colon N \to L$ and $h_1,h_2 \colon L \to \Oh_{amb}$  are four $\Oh_{amb}$-module homomorphisms
with the property
\begin{equation}   \label {eq!assumption}
   h_1 \circ u_2 + h_2 \circ u_1  = 0,
\end{equation}
cf.~(\ref{eq!commutativity2}). We will show by induction on $p+q$ that 
(\ref{eq!basicformula}) holds. If $p = q= 1$ the result is clear, since it
is exactly our hypothesis (\ref{eq!assumption}).

Assume now that $p+q \geq 3$ and that  (\ref{eq!basicformula}) holds
for the $(p-1,q)$ and $(p, q-1)$ cases.

Using Lemma~\ref{lem!firstgoodidentity} we have
\begin {eqnarray*}
    d^{p+q-1}h_2  \circ \alt(u_1^p,u_2^{q-1}) + d^{p+q-1}h_1 \circ
\alt(u_1^{p-1},u_2^{q})  =  \quad \quad \quad \quad \quad \\
     d^{p+q-1}h_2 \circ  [u_1 \wedge  \alt(u_1^{p-1}, u_2^{q-1}) +
          u_2 \wedge  \alt(u_1^p, u_2^{q-2})]   \\ 
    + d^{p+q-1}h_1 \circ [u_1 \wedge  \alt(u_1^{p-2}, u_2^{q}) +
          u_2 \wedge  \alt(u_1^{p-1}, u_2^{q-1})]   
\end {eqnarray*}
which, using (\ref{eqn!onemore}), is equal to 
\begin {eqnarray*}
    (h_2 \circ u_1) \wedge  \alt(u_1^{p-1},u_2^{q-1 })   
       - u_1 \wedge [d^{p+q-2}h_2 \circ \alt(u_1^{p-1},u_2^{q-1})] \\
   + (h_2 \circ u_2) \wedge  \alt(u_1^{p},u_2^{q-2})  
       - u_2 \wedge [d^{p+q-2}h_2 \circ \alt(u_1^{p},u_2^{q-2})]  \\
   + (h_1 \circ u_1) \wedge  \alt(u_1^{p-2},u_2^{q })   
       - u_1 \wedge [d^{p+q-2}h_1 \circ \alt(u_1^{p-2},u_2^{q})] \\
   + (h_1 \circ u_2) \wedge  \alt(u_1^{p-1},u_2^{q-1})  
       - u_2 \wedge [d^{p+q-2}h_1 \circ \alt(u_1^{p-1},u_2^{q-1})]  
\end {eqnarray*}
which, using (\ref{eq!assumption}), is equal to
\begin {eqnarray*}
   -u_1 \wedge[d^{p+q-2}h_2 \circ \alt(u_1^{p-1},u_2^{q-1}) +  
          d^{p+q-2}h_1 \circ \alt(u_1^{p-2},u_2^{q}) ] \\
   -u_2 \wedge [d^{p+q-2}h_2 \circ \alt(u_1^{p},u_2^{q-2}) + 
                d^{p+q-2}h_1 \circ \alt(u_1^{p-1},u_2^{q-1})]  \\
     + (h_2 \circ u_2) \wedge  \alt(u_1^{p},u_2^{q-2}) 
     + (h_1 \circ u_1) \wedge  \alt(u_1^{p-2},u_2^{q })
\end {eqnarray*}
which, using the inductive hypothesis, is equal to 
\begin {eqnarray*}
  (h_2 \circ u_2) \wedge  \alt(u_1^{p},u_2^{q-2}) 
    - u_1 \wedge[\alt(u_1^{p-1},u_2^{q-2})  \circ d^{p+q-2}(h_2 \circ u_2)] \\
    - u_2 \wedge[\alt(u_1^{p},u_2^{q-3})  \circ d^{p+q-2}(h_2 \circ u_2)]  \\
    + (h_1 \circ u_1) \wedge  \alt(u_1^{p-2},u_2^{q}) 
    - u_1 \wedge[\alt(u_1^{p-3},u_2^{q})  \circ d^{p+q-2}(h_1 \circ u_1)] \\
    - u_2 \wedge[\alt(u_1^{p-2},u_2^{q-1})  \circ d^{p+q-2}(h_1 \circ u_1)]
\end {eqnarray*}
which, using  Lemma~\ref{lem!secondgoodidentity}, is equal to
\[
    \alt(u_1^{p-2},u_2^q) \circ
   d^{p+q-1}(h_1 \circ u_1)  + \alt(u_1^{p},u_2^{q-2}) \circ
   d^{p+q-1}(h_2 \circ u_2)
\]
which finishes the proof of Proposition~\ref{prop!basicformula}.

\subsection {The connecting homomorphisms} \label{subs!connhomom}

For an integer $p$, with $1 \leq p \leq n-1$, we set
\[
  T_p = \begin {pmatrix}
   \alt (u_1^p) + z \alt(u_1^{p-2},u_2^2) 
       + z^2  \alt(u_1^{p-4},u_2^4) + \dots  \\
   \alt (u_1^{p-1},u_2) + z \alt(u_1^{p-3},u_2^3)
       + z^2  \alt(u_1^{p-5},u_2^5) + \dots
  \end {pmatrix}, 
\]
with the summation as long as all exponents are nonnegative. More precisely, 
using the notational convention $z^0 = 1$,  and denoting, for $j=1,2$,
by $T_p^j$ the $j$-th row of $T_p$, we have 
\[
  T_p^1=  \underset{i=0}{\overset{[\frac{p}{2}]}{\sum}}  z^i \alt(u_1^{p-2i},u_2^{2i}),
     \quad \quad 
  T_p^2= \underset{i=0}{\overset{[\frac{p-1}{2}]}{\sum}}  z^i \alt(u_1^{p-2i-1},u_2^{2i+1}).
\]

We consider $T_p$ as an $\Oh_{amb}$-homomorphism $\wedge^p N \to \wedge^p L \oplus  \wedge^p L$.
Taking the two rows of $T_p$ we get two $\Oh_{amb}$-homomorphisms  
$\; T_p^1,T_p^2 \colon \wedge^p N \to \wedge^p L$.  For example, 
\begin {eqnarray*}
   T_1 & = &  \begin {pmatrix}
           u_1 \\
           u_2 
     \end {pmatrix}, \  \  
 T_2 =  \begin {pmatrix}
           \wedge^2 u_1 + z \wedge^2 u_2 \\
            u_1 \wedge u_2 + u_2 \wedge u_1 
     \end {pmatrix},  \ \
 T_3 = \begin {pmatrix}  
       \wedge^3 u_1 + z \alt (u_1, u_2^2) \\
     \alt (u_1^2, u_2) + z \wedge^3 u_2    
     \end{pmatrix},       \\
   T_4 & = &\begin {pmatrix}
       \wedge^4 u_1 + z \alt (u_1^2, u_2^2) + 
               z^2  \wedge^4 u_2  \\
     \alt (u_1^3, u_2) + z  \alt (u_1, u_2^3) 
     \end{pmatrix},  \\    
   T_5 & = &\begin {pmatrix}
       \wedge^5 u_1 + z \alt (u_1^3, u_2^2) + 
               z^2   \alt (u_1, u_2^4) \\ 
     \alt (u_1^4, u_2) + z  \alt (u_1^2, u_2^3)
  +z^2 \wedge^5u_2
     \end{pmatrix}.
\end{eqnarray*}

For the following proposition recall that the map
$\phi_p$ was defined in Subsection~\ref{subs!secondcmplx}
while the map $\psi$ was defined in Subsection~\ref{subs!firstcmplx}.

\begin{prop}   \label {prop!bigdiagram}
For any $p$, with $2 \leq p \leq n-1$,   we have a commutative diagram
\begin{equation*}
\begin{CD}
               \wedge^{p-1} N  @< d^p\psi << \wedge^p N \\
       @V T_{p-1} VV     @V  T_p  VV   \\
 \wedge^{p-1} L \oplus  \wedge^{p-1} L  @< \phi_p <<  
      \wedge^{p} L \oplus  \wedge^{p}  L
\end{CD}
\end{equation*}
That is,
 \[
    \phi_p \circ T_p = T_{p-1} \circ d^p\psi
 \]
 as maps  $\; \; \; \wedge^p N \to \wedge^{p-1} L \oplus \wedge^{p-1} L$.
\end{prop}

\begin {pf}
   Fix $p \geq 2$.  For $i \in \Z$ we set  
\[
        C_i = \alt (u_1^i, u_2^{p -i})
\]
if $ 0 \leq i \leq p$ and $C_i =0$  otherwise, and we set
\[
        D_i = \alt (u_1^i, u_2^{p-1 -i})
\]
if $ 0 \leq i \leq p-1$ and $ D_i =0$ otherwise.  

For the rest of the proof all sums are for $i \in \Z$.

We have 
\begin {eqnarray*}
       T_p^1  =    \sum z^i C_{p-2i}, \quad \quad  \quad T_p^2  =   \sum z^i C_{p-2i-1}, \\
       T_{p-1}^1  =   \sum z^i D_{(p-1)-2i},	  \quad \quad \quad T_{p-1}^2  =    \sum z^i D_{(p-1)-2i-1}.
\end {eqnarray*}
By Proposition~\ref{prop!basicformula}, for any $i \in \Z$ we
have
\begin {equation}  \label{basiceqnforCD}
    d^ph_1 \circ C_{i-1} + d^ph_2 \circ C_{i} = D_{i-2} \circ d^p(h_1 \circ u_1) +   
                D_{i} \circ d^p(h_2 \circ u_	2).
\end {equation}
Using (\ref{eq!commutativity1}), we need to show  that
\[
    d^ph_1 \circ T_p^1 + z d^ph_2 \circ T_p^2  = T_{p-1}^1 \circ d^p(h_1 \circ  u_1+zh_2 \circ u_2) 
\]
and that
\[
    d^ph_2 \circ T_p^1 + d^ph_1 \circ T_p^2 = T_{p-1}^2 \circ d^p(h_1 \circ u_1+zh_2 \circ u_2).
\]
Using (\ref{basiceqnforCD}) we get
\[
     d^ph_1 T_p^1 + z d^ph_2 T_p^2 = \sum [z^i (d^ph_1 C_{p-2i} + d^p h_2 C_{p-2i+1})] 
\]
\[  
    =  \sum [z^i (D_{p-2i-1} d^p(h_1u_1)  + D_{p-2i+1} d^p(h_2u_2))]
\]
\[
   =   \sum [z^i (D_{(p-1)-2i} (d^p(h_1u_1)+ zd^p(h_2u_2)))] = T_{p-1}^1 d^p(h_1u_1+zh_2u_2).
\]

Similarly,
\[
    d^ph_2 T_p^1 + d^ph_1 T_p^2 =  \sum [z^i (d^ph_2 C_{p-2i} + d^ph_1 C_{p-2i-1})] 
\]
\[
  =  \sum [z^i (D_{p-2i-2} d^p(h_1u_1)  + D_{p-2i} d^p(h_2u_2)) ]
\]
\[
  =  \sum [z^i (D_{(p-1)-2i-1} (d^p(h_1u_1)+ zd^p(h_2u_2)))] = T_{p-1}^2 d^p(h_1u_1+zh_2u_2).
\]
\medskip \QED \end {pf}

\subsection {Linear equations of type $\II_1$ unprojection}

Assume $n \geq 2$ is given. We use the notations of the previous subsections.
In Subsection \ref{subs!connhomom} we defined two $\Oh_{amb}$-module homomorphisms  
$\quad T_{n-1}^1, T_{n-1}^2 \colon $ \newline $\wedge^{n-1}N \to \wedge^{n-1}L$.

There exist, for $j=1,2$ and $1 \leq i \leq n$, unique element $\sigma_i^j \in \Oh_{amb}$ such that
\[ 
  T_{n-1}^j  (\widetilde{e}_1 \wedge \dots \wedge \widetilde{e}_{n-1})
   = \sum_{i=1}^n (-1)^{i-1} \sigma_i^j  e_1 \wedge \dots \wedge \widehat{e_i}
    \wedge \dots      \wedge e_n.
\]  
For $1 \leq i \leq n$ we define the polynomials $l_i, l_{i+n} \in \Oh_{amb}[s_0,s_1]$ with 
\begin{eqnarray}
   l_i     =  zx_i s_0 + y_i s_1 +\sigma_i^1      \notag  \\ 
   l_{i+n} =  y_i s_0 + x_i s_1 + \sigma_i^2      \label {eqn!linearsecond}
\end{eqnarray}

Using Proposition~\ref{prop!bigdiagram} and arguing as in \cite{R}~Section~9.5 
we get the following Theorem.

\begin {theorem}   \label{thm!aboutlinears} Fix a parameter value $n \geq 2$.
The polynomials $l_1, \dots ,l_{2n}  \in \Oh_{amb}[s_0,s_1]$ 
specified in (\ref{eqn!linearsecond}) are the linear relations of the 
type $\II_1$ unprojection in the sense  of   Remark~\ref{rem!aboutIY}.
\end {theorem}

\section {Quadratic relation of Type $\II_1$ for $n=3$}  \label{sec!quadratic}

Assume now that we are in the type $\II_1$ unprojection with $n=3$.
We follow the notations of  Section~\ref{sec!notation}
with only one change: 
Since our symmetric formula (\ref{eqn!dfnofb}) 
for the quadratic relation will need $1/2$ as
coefficient, our ambient ring $\Oh_{amb}$ will now be 
\[
      \Oh_{amb} = \Z[\frac {1}{2}] [x_1, \dots , x_3, y_1, \dots ,y_3,z, A^p_{ij},B^p_{lm} ],
\]
with indices $1 \leq p \leq 2$ and $1 \leq i < j \leq 3, \; 1 \leq l \leq  m \leq 3$.

Recall that in   Section~\ref{sec!notation} we defined  two symmetric 
$3 \times 3$ matrices $B^1$ and  $B^2$ by
\[
   B^j = 
        \begin {pmatrix}
             B_{11}^j & B_{12}^j & B_{13}^j   \\ 
                      & B_{22}^j & B_{23}^j   \\
            \text {sym}  &  & B_{33}^j  
   \end{pmatrix},  
\] 
for $j=1,2$, two skew--symmetric $3 \times 3$ matrices 
$A^1$ and  $A^2$ by
\[
   A^j = 
        \begin {pmatrix}
                 0    &  A_{12}^j & A_{13}^j   \\ 
                      &  0        & A_{23}^j   \\
               \text {--sym}    &   &   0   
   \end{pmatrix},  
\] 
for $j=1,2$, and two $1 \times 3$ matrices $x,y$ with
\[
   x = (x_1,x_2,x_3), \ \  y = (y_1,y_2,y_3).
\]

We now define three induced $3 \times 3$ symmetric matrices $\ad B^1,  \;
\ad B^2$ and $\ad B^{1,2}$.

For a square $n \times n$ matrix $M$ with entries
in a commutative ring $R$ with unit, we  denote by $\ad M$ the 
 $n \times n$ matrix  
with $kl$-th entry, for $1 \leq k,l \leq n$,  equal to  $(-1)^{k+l}$ times
the determinant of the submatrix of $M$ obtained by deleting the 
$l$-th row and the $k$-th column. It is well-known that we have   
the identities  $M (\ad M) = (\ad M) M = (\det M)I_n$,
where $I_n$ is the identity $n \times  n$ matrix over $R$.

In our situation, since $B^1$ and $B^2$ are symmetric $3 \times 3$  
matrices we have that
$\ad B^1$ and $\ad B^2$  are symmetric $3 \times 3$ matrices, such that
for $j=1,2$
\[
  \ad B^j =  \begin {pmatrix}
            B^{j}_{22}B^{j}_{33}-(B^{j}_{23})^2
        & -(B^{j}_{12}B^{j}_{33}-B^{j}_{13}B^{j}_{23})
        &   B^{j}_{12}B^{j}_{23}-B^{j}_{13}B^{j}_{22} \\
          
        &   B^{j}_{11}B^{j}_{33}-(B^{j}_{13})^2
        & -(B^{j}_{11}B^{j}_{23}- B^{j}_{13}B^{j}_{12} )  \\
        \text {sym}   
        &   
        &  B^{j}_{11}B^{j}_{22} -(B^{j}_{12})^2
          \end {pmatrix}.
\]   
Notice that for $1 \leq k \leq 3,  \;
1 \leq l \leq 3$, we have
\[
    (\ad B^1)_{kl}  = (-1)^{k+l} ( B^1_{ij}B^1_{mn}  -B^1_{in}B^1_{mj})
\]
where $\{ i,k,m \} = \{ j,l,n \} = \{1,2,3\}$ and $i<m,\ j<n$, and
similarly for $\ad B^2$.

We define  $\ad B^{1,2}$ to be the $3 \times 3$ symmetric matrix, such that,
for $1 \leq k \leq 3,  \;1 \leq l \leq 3$, 
\[
    (\ad B^{1,2})_{kl}  = (-1)^{k+l} ( B^1_{ij}B^2_{mn} +  B^2_{ij}B^1_{mn} -
                      B^1_{in}B^2_{mj} -  B^2_{in}B^1_{mj}),
\]
where $\{ i,k,m \} = \{ j,l,n \} = \{1,2,3\}$ and $i<m,\ j<n$.
By the above definition we have
\begin {eqnarray*} 
   \ad B^{1,2}_{11}  & = &  B^1_{22}B^2_{33}+B^2_{22}B^1_{33}-2B^1_{23}B^2_{23}\\
   \ad B^{1,2}_{12}  & = & -B^1_{12}B^2_{33}-B^2_{12}B^1_{33}+B^1_{13}B^2_{23}+B^2_{13} B^1_{23}\\
   \ad B^{1,2}_{13}  & = &  B^1_{12}B^2_{23}+ B^2_{12}B^1_{23}-B^1_{13}B^2_{22}-B^2_{13}B^1_{22} \\
   \ad B^{1,2}_{22}  & = &  B^1_{11}B^2_{33}+B^2_{11}B^1_{33}-2 B^1_{13}B^2_{13} \\
   \ad B^{1,2}_{23}  & = & -B^1_{11}B^2_{23}-B^2_{11}B^1_{23}+B^1_{13}B^2_{12}+B^2_{13} B^1_{12}\\
   \ad B^{1,2}_{33}  & = &  B^1_{11}B^2_{22}+B^2_{11}B^1_{22}-2B^1_{12}B^2_{12}
\end {eqnarray*}
and the rest of the entries of $\ad B^{1,2}$ are determined by symmetry.

We now define two polynomials $a,b \in \Oh_{amb}$.
We set 
\begin{eqnarray}  \label{eqn!dfnofa}
     a &=&  \det  \begin {pmatrix}
                 x_1    &  x_2 &  x_3  \\ 
               A_{23}^1 &  -A_{13}^1 &  A_{12}^1   \\
               A_{23}^2 &  -A_{13}^2 &  A_{12}^2   
   \end{pmatrix}  \\ \notag
        & = & x_1 (A_{12}^1A_{13}^2-A_{12}^2A_{13}^1)
         +x_2 (A_{12}^1A_{23}^2-A_{12}^2A_{23}^1) \\ \notag
         &  & + \phantom{=} x_3 (A_{13}^1A_{23}^2-A_{13}^2A_{23}^1), 
\end{eqnarray}
and
\begin{eqnarray}  \label{eqn!dfnofb}
     b &=&  \frac{1}{2} \Bigl[  yB^1 (\ad B^2)B^1 y^t +  yB^2 (\ad B^1) B^2y^t \\ \notag
       &  & \phantom{=}   +  \phantom{=}  2 xA^1 (\ad B^2) (A^1)^t x^t +2 xA^2 (\ad B^1) (A^2)^t x^t
               \\ \notag
       &  &  \phantom{=} -\phantom{=}  2  x A^1 (\ad B^{1,2})(A^2)^tx^t \\ \notag
       &  & \phantom{=}  - \phantom{=} z  xB^1 (\ad B^2)B^1 x^t -  zxB^2 (\ad B^1) B^2x^t \\ \notag
       &  &\phantom{=}  +  \phantom{=} 4  yB^1 (\ad B^2)(A^1)^t x^t + 4 yB^2 (\ad B^1) (A^2)^t x^t
                 \\ \notag
       &  &\phantom{=}  -\phantom{=}  ( \sum_{1\leq i<j \leq 3}A_{ij}^1 (x_iy_j-x_jy_i) )
                (\sum_{1 \leq i,j \leq 3}^{3} B_{ij}^1 \ad B_{ij}^2)   \\ \notag
       &  & \phantom{=}  - \phantom{=} ( \sum_{1\leq i<j \leq 3}A_{ij}^2 (x_iy_j-x_jy_i) )
               ( \sum_{1\leq i,j \leq 3} B_{ij}^2 \ad B_{ij}^1 )  \Bigr]. 
\end{eqnarray}

The two polynomials $f_1,f_2$ corresponding to the ones in (\ref{eq!dfnoff}), are
\begin{equation}  
   f_j = xA^jy^t + yB^jy^t - zxB^jx^t, \label {eqn!dfnofp1p2}
\end{equation}
for $j=1,2$.

We denote by $u_1$  the $3 \times 2$ matrix
\[
   u_1 = (-A^1x^t + B^1y^t  \; |  -A^2x^t  + B^2y^t)
\]
and by $u_2$  the $3 \times 2$ matrix
\[
  u_2 = (-B^1x^t \; | -B^2x^t).
\]
These matrices
correspond to the maps with the same name defined in~(\ref{eqn!matrixformat}).

For $j =1,2$, we set  $\widetilde{\phi_1}(u^j)$ to be the $3 \times 1$ matrix
\[
    \widetilde{\phi_1}(u^j)  =  \begin {pmatrix}  
                  u^j_{21}u^j_{32}-u^j_{22}u^j_{31} \\
                 -u^j_{11}u^j_{32}+u^j_{12}u^j_{31} \\ 
                  u^j_{11}u^j_{22}-u^j_{12}u^j_{21}
               \end {pmatrix}
\]
and we also define the $3 \times 1$ matrix
\[
   \widetilde{\phi_2}(u^1,u^2) = \begin {pmatrix}  
                      u^1_{21}u^2_{32}+u^2_{21}u^1_{32}-u^1_{22}u^2_{31}-u^2_{22}u^1_{31}  \\
                     -u^1_{11}u^2_{32}-u^2_{11}u^1_{32}+u^1_{12}u^2_{31}+u^2_{12}u^1_{31} \\
                      u^1_{11}u^2_{22}+u^2_{11}u^1_{22}-u^1_{12}u^2_{21}-u^2_{12}u^1_{21}
               \end {pmatrix}.
\]
The matrix $\widetilde{\phi_1}(u^j)$ corresponds to the map $ \wedge^2 u_j$ and the matrix 
$\widetilde{\phi_2}(u^1,u^2)$  corresponds to the map $\alt(u_1,u_2)$ which were defined
in Subsection~\ref{sub!basicformula}.

We denote, for $1 \leq i \leq 3$, 
\begin {eqnarray}
     l_i &  = &  zx_i s_0 + y_i s_1 + [\widetilde{\phi_1}(u^1)+z\widetilde{\phi_1}(u^2)]_{i1}  \label{eqn!dfnoflinears}  \\
     l_{i+3} & = &  y_i s_0 + x_i s_1 + [\widetilde{\phi_2}(u^1,u^2)]_{i1},  \notag
\end{eqnarray}
where we use the usual notation $M_{ij}$ for the $ij$-th entry of a matrix $M$. 
Clearly, for $1 \leq i \leq 6$, $\; l_i \in \Oh_{amb}[s_0,s_1]$ is affine linear with respect
to $s_0$ and $s_1$.
The  polynomials  $l_1, \dots ,l_6$ correspond to the polynomials with the same name in (\ref{eqn!linearsecond}),
hence by Theorem~\ref{thm!aboutlinears} they are the linear relations of the unprojection in
the sense of   Remark~\ref{rem!aboutIY}.

The main result of this section is the following theorem, which 
specifies the quadratic relation of the unprojection ring for the
parameter value $n=3$. It 
will be proved in Subsection~\ref{sub!pfofquadratic}. 

\begin {theorem}  \label{thm!quadratic}
The polynomial 
\begin {equation}  \label{eqn!dfnofq}
     q  =    s_1^2 - z s_0^2 - a s_0 + b \in \Oh_{amb}[s_0,s_1],
\end {equation} 
where $a$ was defined in (\ref{eqn!dfnofa}) and  $b$ was defined 
in (\ref{eqn!dfnofb}),
is the quadratic  relation of the type $\II_1$ unprojection for the parameter
value $n=3$ in the sense of  Remark~\ref{rem!aboutIY}.

In other words,  the ideal  $I_Y \subset \Oh_{amb}[s_0,s_1]$ defining,
in the sense of Section~\ref{sec!notation}, the type $\II_{1}$  unprojection 
is equal to 
\[    
            I_Y = (f_1,f_2) + (l_1, \dots , l_6) + (q),
\]
where the two polynomials $f_1,f_2$ generating $I_X$ were  defined 
in (\ref{eqn!dfnofp1p2}) and  the six polynomials $l_1,  \dots ,l_6$ were  
calculated  in (\ref{eqn!dfnoflinears}).
\end{theorem}

\subsection { Proof of Theorem~\ref{thm!quadratic} }  \label{sub!pfofquadratic}

We first prove the following reduction lemma.

\begin {lemma}   \label {lem!reduction}
  To prove Theorem \ref{thm!quadratic} it is enough to prove that
that the element $(x_1+x_2+x_3) q$ of $\Oh_{amb}[s_0,s_1]$ is inside the ideal 
 $(f_1,f_2,l_1, \dots ,l_6)$ of $\Oh_{amb}[s_0,s_1]$.
\end {lemma}

\begin {pf}
   Indeed, by \cite{P}~Proposition~2.16  there exists  
$\widetilde {q} \in  \Oh_{amb}[s_0,s_1]$ of the form
\[
   \widetilde {q} = s_1^2-zs_0^2 + \widetilde{u}
\]
where $\widetilde{u}$ is affine linear in $s_0,s_1$,  such that the ideal of the 
unprojection ring is equal to  $J$, with
\[
   J = (f_1,f_2,l_1, \dots ,l_6, \widetilde{q}) \subset  \Oh_{amb}[s_0,s_1].
\]
Since by the assumptions of the lemma $(x_1+x_2+x_3) q  \in J$, 
and by~\cite{P}~Section~3.1 $J$ is a prime ideal of $\Oh_{amb}[s_0,s_1]$ 
we get $q \in J$. Therefore, the affine linear with respect to $s_0,s_1$ 
element $q -\widetilde{q}$ is in $J$, as a consequence
(compare proof of \cite{P}~Proposition~2.13) we get
\[ 
     q- \widetilde {q} \in (f_1,f_2,l_1, \dots ,l_6),
\]
hence 
\[  
   J = (f_1,f_2,l_1, \dots ,l_6, q)
\]
which finishes the proof of Lemma~\ref{lem!reduction}.
\medskip \QED  \end{pf}

For the rest of the proof we define some more notation:

We define two polynomials $\widetilde{g_1},\widetilde{g_2}$, such that, for $j \in \{1,2 \}$, we have
\begin {eqnarray*}  
 \widetilde{g_j} & =  & \phantom {+} (x_1  B^j_{11}+ x_2 B^j_{12} + x_3  B^j_{13}) (\ad B^p_{11} + \ad B^p_{12} +\ad B^p_{13}) \\
    &  &  +   (x_1  B^j_{21}+ x_2  B^j_{22} + x_3  B^j_{23})  (\ad B^p_{21} + \ad B^p_{22} +\ad B^p_{23})    \notag  \\ 
    &  &  +   (x_1  B^j_{31}+ x_2  B^j_{32} + x_3  B^j_{33})   (\ad B^p_{31} + \ad B^p_{32} +\ad B^p_{33})     \notag 
\end {eqnarray*}
where $p \in \{1,2 \}$ is uniquely specified by $\{j,p \} = \{1,2 \}$.

Moreover, we define two $1 \times 3$ matrices 
\[
   L_1 = (l_1,l_2,l_3), \quad    L_2 = (l_4,l_5,l_6)
\]
and two $3 \times 3$ matrices $M_1$ and $M_2$ by
\[
  (M_1)_{ij} = \widetilde{\phi_3} (B^1_{i1}, B^1_{i2}, B^1_{i3}, B^1_{j1}, B^1_{j2},B^1_{j3},
                     B^2_{i1}, B^2_{i2}, B^2_{i3}, B^2_{j1}, B^2_{j2},B^2_{j3})
\]
and
\begin {eqnarray*}
  (M_2)_{ij} & = & \phantom {-} \widetilde{\phi_3} (A^1_{i1}, A^1_{i2}, A^1_{i3}, A^1_{i1}, A^1_{i2},A^1_{i3},
                     B^2_{j1}, B^2_{j2}, B^2_{j3}, B^2_{j1}, B^2_{j2},B^2_{j3}) \\
                  & &     -\widetilde{\phi_3} (A^2_{i1}, A^2_{i2}, A^2_{i3}, A^2_{i1}, A^2_{i2},A^2_{i3},
                     B^1_{j1}, B^1_{j2}, B^1_{j3}, B^1_{j1}, B^1_{j2},B^1_{j3})
\end {eqnarray*}
where
\begin {eqnarray*}
   \widetilde{\phi_3} (a_1,a_2,a_3,b_1,b_2,b_3,  c_1,c_2,c_3, d_1,d_2,d_3) =    (a_1(d_3-d_2) + a_2  (d_1-d_3)    \\ 
      + a_3  (d_2-d_1))  -(c_1(b_3-b_2) + c_2  (b_1-b_3) + c_3  (b_2-b_1)).
\end {eqnarray*}
Finally, we define the polynomial 
\begin {eqnarray*}  
   \widetilde{w} = x M_1 (L_1)^t + y M_1 (L_2)^t +  x M_2 (L_2)^t.
\end {eqnarray*}

\begin {prop} \label{prop!mainidentity}
We have the following equality
\begin {equation}   \label{eqn!bigcheck}
   2(x_1+x_2+x_3)q = -\widetilde{w}+\widetilde{g_1}f_1+\widetilde{g_2}f_2 - 2s_0 (l_1+l_2+l_3) + 2s_1(l_4+l_5+l_6)
\end {equation} 
of elements of  $\Oh_{amb}[s_0,s_1]$.
\end {prop}

\begin {pf}
The verification was done using the computer algebra program Macaulay 2
\cite{GS93-08}.
\medskip \QED  \end{pf}

The proof of Theorem~\ref{thm!quadratic} follows now by combining 
Proposition~\ref{prop!mainidentity} with  Lemma~\ref{lem!reduction}.

\subsection {Invariance under $\SL_3$}  \label{subs!invariance}

Let $P \in \SL_3(\Z)$ be a $3 \times 3$ matrix with integer coefficients and
determinant one.  Define  two $1 \times 3$ matrices 
    $\widetilde {x} = (\widetilde{x_1},\widetilde{x_2},\widetilde{x_3})$ and 
    $\widetilde {y} = (\widetilde{y_1},\widetilde{y_2},\widetilde{y_3})$
by
\[
      \widetilde {x} = xP^{-t}, \quad \quad    \widetilde {y} = yP^{-t},
\]
(where $P^{-t}\;$ is, of course, the transpose matrix of the inverse of $P$),
two $3 \times 3$ skew--symmetric matrices    $\widetilde {A^j}$ 
and two $3 \times 3$ symmetric matrices    $\widetilde {B^j}$,
for $j=1,2$, by 
\[
      \widetilde {A^j} = P^tA^jP,  \quad \quad   \widetilde {B^j} = P^tB^jP
\]
and finally set $\widetilde {z} = z$.
Moreover, we denote by $\widetilde  {A^j_{lm}}$ the $lm$-th entry of the matrix
$\widetilde {A^j}$  and by $\widetilde  {B^j_{lm}}$ the $lm$-th entry of the matrix
$\widetilde {B^j}$.

It can be checked that if we denote by $\widetilde{a}$ and  $\widetilde{b}$ the polynomials
of Expressions~(\ref{eqn!dfnofa}) and (\ref{eqn!dfnofb}) respectively,  with the variables 
without tilde replaced by the corresponding expressions with tilde,   we have  
$ \widetilde{a} = a$ and $\widetilde{b} = b$.

\subsection { How we got to Theorem~\ref{thm!quadratic} }  \label{sub!howwegot}

In this subsection we briefly sketch the method that allowed us
to calculate, using the computer program Maple, the expressions for $a$ and $b$ 
appearing in Theorem~\ref{thm!quadratic}. 

{\it  Step 1.}  Recall the explicit expressions for $l_i$ in  (\ref{eqn!dfnoflinears}).
 We started from the expression 
\[ 
    2\;[-s_0(l_1+l_2+l_3)  + s_1 (l_4+l_5+l_6)]  
\]
which is equal to 
\[
   2 (x_1+x_2+x_3) (s_1^2-zs_0^2) + 2 H,
\]
where
\[
    H =  -s_0\sum_{i=1}^3 [\phi_1(u^1)+z\phi_1(u^2)]_{i1} +
                  s_1 \sum_{i=1}^3  [\phi_2(u^1,u^2)]_{i1}.
\]
It is easy to see that $2H$ is a homogeneous quadratic
polynomial with respect to the variables $A^t_{ij}, B^s_{kl}$.
More precisely,  $2H$  has a (unique) natural representation as 
a sum of terms of the form  $ \quad \quad    A^1_{ij} A^2_{kl} c^A_{ijkl} \quad \quad $ 
plus sum of terms of the form $ \quad  \quad   A^1_{ij} B^2_{kl} c^{AB}_{ijkl}  \quad \quad  $
plus sum of terms of the form $ \quad  \quad    A^2_{ij} B^1_{kl} c^{BA}_{ijkl} \quad  \quad$
plus sum of terms of the form  $ \quad  \quad    B^1_{ij} B^2_{kl} c^{B}_{ijkl}$.

{\it Step 2.}    We first handle the 'subpart' of $2H$ which is sum 
of terms of the form  $   A^1_{ij} A^2_{kl} c^A_{ijkl} $.
These term appear only as contributions from  $s_0 (l_1+l_2+l_3)$
and, moreover, this subpart is equal
to $2as_0$, where $a$  was defined in (\ref{eqn!dfnofa}).

{\it Step 3.}   We now handle the 'subpart' of $2H$ which is sum 
of terms of the form  $   A^1_{ij} B^2_{kl} c^{AB}_{ijkl} $. 
We notice that there is a kind of complementarity of terms
for the pairs $(l_1,l_4), (l_2,l_5), (l_3,l_6)$ in the sense
of the following example:  
  
  A partial sum of the coefficient 
of $A^1_{13}B^2_{12}$ coming from $l_1$ and $l_4$ 
is $2y_1x_1s_0 +2x_1^2s_1$ which can be written as  
$2x_1 (y_1s_0+x_1s_1)$. Using $l_4$ we get an 
expression not involving $s_0$ and $s_1$. Similarly,
another partial sum  coming  from $l_2$ and $l_5$ is
$ 2x_1(y_2s_0+x_2s_1)$ and we can use $l_5$ to 
to get an expression not involving $s_0$ and $s_1$.
A similar procedure can be done for the the partial sum
coming from $l_3$ and $l_6$.

By symmetry, the handling of the 'subpart' of $2H$ which is sum 
of terms of the form $A^2_{ij} B^1_{kl} c^{BA}_{ijkl}$ is similar.

The handling of the 'subpart' of $2H$ which is sum 
of terms of the form   $B^1_{ij} B^2_{kl} c^{B}_{ijkl}$ is similar,
but we need to take care of preserving the symmetries.
This is the only part where we actually need the coefficient $2$ in $2H$.

To give an example,  a partial sum of the coefficient 
of $B^1_{13} B^2_{23}$  coming from $l_1$ and $l_4$ 
is 
\[
    2\;[-s_0(-x_1x_3z-y_1y_3)+s_1(y_3x_1+y_1x_3)]
\]
which we first write as 
\[     
   y_3  (y_1s_0+x_1s_1) +  x_1  (zx_3s_0+y_3s_1) +  
                       y_1  (y_3s_0+x_3s_1)  + x_3  (zx_1s_0 + y_1 s_1)
\]                  
and then we use the appropriate linear equations $l_i$ to get
an expression  not involving $s_0$ and $s_1$.

{\it Step 4.} At the end of step 3 we arrived to an expression not involving
$s_0,s_1$. By suitable subtraction (in a symmetric way)
of multiples of $f_1$ and $f_2$ we get 
a polynomial divisible by $x_1+x_2+x_3$. The quotient is
$b$, which after some further term by term effort can be written
in the form (\ref{eqn!dfnofb}).

\subsection {Type $\II_1$ unprojection for $n=2$} \label{sub!reidcasenis2}

This subsection  contains  the explicit equations obtained by Miles Reid 
for the type $\II_1$ unprojection 
with  parameter value $n=2$. It is taken from \cite{R} Section~9.5.

For $n=2$ the ambient ring is 
\[
    \Oh_{amb}= \Z [x_1,x_2, y_1, y_2, z,  A_{11}, B_{11}, B_{12}, B_{22}],
\]
the matrix $M$ corresponding to the one in (\ref{eq!dfnofM}) is equal to 
\begin{equation*}     
    M =  \begin {pmatrix}
             y_1 & y_2  &  z x_1 & z x_2  \\
             x_1 & x_2  &  y_1   & y_2
     \end {pmatrix},
\end{equation*} 
there is a unique polynomial 
\[
   f=A_{12}(x_1y_2-x_2y_1)+B_{11}(y_1^2-zx_1^2)+2B_{12}(y_1y_2-zx_1x_2)+
      B_{22}(y_2^2-zx_2^2)
\]
corresponding to (\ref{eq!dfnoff}),  the linear equations  $l_1, \dots ,l_4$
are given by
\begin {eqnarray*} 
     l_1 & = &   zx_1s_0 + y_1s_1 +( x_1A_{12}+y_1B_{12}+y_2B_{22})    \\
     l_2 & = &   zx_2s_0 + y_2s_1 +( x_2A_{12}-y_1B_{11}-y_2B_{12})    \\
     l_3 & = &   y_1s_0+  x_1s_1 + ( -x_1B_{12}-x_2B_{22} )          \\
     l_4 & = &   y_2s_0 +  x_2s_1  + (x_1B_{11}+x_2B_{12})
\end {eqnarray*}
and the quadratic equation $q$ is given by
\[
  q  =  s_1^2-zs_0^2-A_{12}s_0-(B_{12})^2+B_{11}B_{22}.
\]

In other words, using the notations of Section~\ref{sec!notation} we
have the equality of ideals of  $\Oh_{amb}[s_0,s_1]$
\[
   I_Y = (f) + (l_1, \dots ,l_4) +  (q).
\]

\begin {rem}
Since  $\quad \quad f = y_2l_1 -y_1l_2+zx_2l_3 -zx_1l_4, \quad \quad$ we even get
\[
    I_Y = (l_1, \dots ,l_4) +  (q).
\]
\end {rem}

\begin {rem} By \cite {P}  the ideal $I_Y$ is Gorenstein codimension $3$. It is easy to
see that it is equal to the ideal generated by the $5$ submaximal Pfaffians of 
the $5 \times 5$  skew--symmetric matrix (with entries in  $\Oh_{amb}[s_0,s_1]$)
\[
   \begin {pmatrix}
        0 & x_1 & x_2 & y_1     &   y_2     \\
          & 0  & -s_0 &  -B_{22} &  s_1+B_{12} \\
          &    &  0  &  -s_1 + B_{12} & -B_{11} \\
          &    &     &        0  & -z s_0- A_{12} \\
          &   \text {--sym}     &     &           &   0
    \end{pmatrix}.
\]
(For a discussion about the Pfaffians of a skew--symmetric matrix see, for example,
\cite{BH}~Section~3.4.)
\end {rem}

\section { Applications to algebraic geometry }   \label{sec!applications}

We believe that the explicit formulas of Theorem~\ref{thm!quadratic}
can be used together with Gavin Brown's online database  
of graded rings \cite{Br} for the proof of the existence of a number of
(singular) Fano $3$-folds with anticanonical ring of codimension $4$.
We discuss below two such examples, the first of which is due to Selma
Alt{\i}nok, and has the interesting
property $| -K_X | = \emptyset$.
We also expect that the explicit formulas of Theorem~\ref{thm!quadratic}
together with the orbifold Riemann--Roch theorem obtained in \cite{BS}
can lead to the construction of new codimension $4$ Calabi--Yau $3$-folds.

\subsection {The example of Alt{\i}nok} \label{subs!altinokexample}

This example is taken from \cite {R} Example~9.14 and is due to Alt{\i}nok.
Consider the weighted projective space $\PP^2(1,3,5)$ with coordinates
$u,v,w$  and the weighted projective space $\PP^5(2,3,4,5,6,7)$ with
coordinates $x,v,y,w,z,t$. We define the map 
$\PP(1,3,5)\to \PP(2,3,4,5,6,7)$  given by
 \[
 x=u^2,\quad v=v,\quad y=uv,\quad w=w,\quad z=uw,\quad t=u^7.
 \]
The equations of the image $D$ of the map are
 \[
 \rank\begin{pmatrix}
 y &  z  &  t    &   xv  &  xw &   x^4 \\
 v &  w  &  x^3  &   y   &  z  &   t
 \end{pmatrix}\leq 1,
 \]
that is,
\[
 \begin{array}{l}
  yw=vz, \quad yx^3=vt, \quad  zx^3=wt,\quad  y^2=v^2x,  \quad  yz=vwx  \\   
  yt=vx^4, \quad z^2=w^2x,  \quad zt=wx^4, \quad t^2=x^7.
 \end{array}
 \]

A general 
complete intersection $X_{12,14}$ containing $D$ is obtained by choosing two general
combinations,   one of degree $12$ and the other of degree $14$,
of the above equations of $D$.

We perform a type $\II_1$ unprojection of the pair $D \subset X_{12,14}$ to
get a codimension $4$ $3$-fold  
\[
    Y \subset  \PP^7(2,3,4,5,6,7,8,9).
\]
Notice that the new variable $s_0$ has degree $8$ and the 
new variable $s_1$ has degree $9$.

After substituting to the formulas (\ref{eqn!linearsecond})
and (\ref{eqn!dfnofq}) we checked the quasismoothness of $Y$ using the
computer algebra program Singular \cite {GPS01}.

\subsection {The second Fano example }

In this subsection we sketch a construction, suggested by Brown's online database  
of graded rings \cite{Br}, of a codimension $4$ Fano $3$-fold 
\[
       Y \subset \PP^7(1,1,2,2,2,2,3,3) 
\]   
starting from a (nongeneric) codimension $2$ complete intersection Fano 
$3$-fold $X_{4,6} \subset  \PP^5(1,1,2,2,2,3)$.

Write  $x_1,x_2,y_1,y_2,y_3,w$ for the coordinates of  $\PP^5(1,1,2,2,2,3)$, 
and consider the subscheme $D \subset \PP^5(1,1,2,2,2,3)$ with equations 
 \[
    \rank\begin{pmatrix}
            y_1    &    y_2    & w     &   y_3x_1  & y_3x_2    & y_3^2     \\
            x_1    &    x_2    & y_3   &   y_1     & y_2       & w
    \end{pmatrix}\leq 1.
\]
By definition, the ideal $I_D$ of $D$ is generated by the $2 \times 2$ minors of the above matrix.
Denote by $X_{4,6} \subset  \PP^5(1,1,2,2,2,3)$ a general codimension
$2$ complete intersection (with equations of degrees  $4$ and $6$) containing $D$. 
The equations of $X_{4,6}$ are obtained by choosing two general
combinations,   one of degree $4$ and the other of degree $6$,
of the above equations of $D$.

We perform a type $\II_1$ unprojection of the pair $D \subset X_{4,6}$ to
get a codimension $4$ $3$-fold  $ Y \subset  \PP^5(1,1,2,2,2,2,3,3)$.
Notice that here the new variable $s_0$ has degree $2$ and the new variable 
$s_1$ has degree $3$.

After substituting to the explicit formulas 
 (\ref{eqn!linearsecond}) and (\ref{eqn!dfnofq})
we checked the quasismoothness of $Y$  using the
computer algebra program Singular \cite {GPS01}.

\section {Appendix: Macaulay 2 code}   \label{sec!appendix}

This section contains the Macaulay 2 code for the type $\II_1$ unprojection
with $n=3$.

\begin{verbatim}
-- M2 code for type II_1 unprojection for n=3

kk = QQ
S = kk [A1p12,A1p13,A1p23,A2p12,A2p13,A2p23, B1p11,B1p12,B1p13,
        B1p22,B1p23,B1p33, B2p11,B2p12,B2p13,B2p22,B2p23,B2p33, 
        y1,y2,y3,x1,x2,x3,z,s0,s1, 
        Degrees => {6:2, 12:1, 3:2, 3:1,2,3,4}] 

--adjointMatrix: input:   M  3x3 matrix, 
--               output:  classical adjoint of M (which is 3x3)
adjointMatrix  =  (M) -> transpose matrix {
  {det submatrix(M,{1,2},{1,2}), -det submatrix(M,{1,2},{0,2}),
              det submatrix(M,{1,2},{0,1})},
  {-det submatrix(M,{0,2},{1,2}), det submatrix(M,{0,2},{0,2}),
             -det submatrix(M,{0,2},{0,1})}, 
  {det submatrix(M,{0,1},{1,2}),-det submatrix(M,{0,1},{0,2}),
              det submatrix(M,{0,1},{0,1})}}

--doubleAdjointMatrix: input:   M1, M2 are  3x3 matrices
--                     output:  3x3 matrix 
doubleAdjointMatrix  =  (M1,M2) -> (
   M111 := M1_0_0; M112 := M1_1_0; M113 := M1_2_0;
   M121 := M1_0_1; M122 := M1_1_1; M123 := M1_2_1;
   M131 := M1_0_2; M132 := M1_1_2; M133 := M1_2_2;
   M211 := M2_0_0; M212 := M2_1_0; M213 := M2_2_0;
   M221 := M2_0_1; M222 := M2_1_1; M223 := M2_2_1;
   M231 := M2_0_2; M232 := M2_1_2; M233 := M2_2_2;
   transpose matrix {  
       {M122*M233-M123*M232+ M222*M133-M223*M132,
       -(M121*M233-M123*M231+ M221*M133-M223*M131 ),
       M121*M232-M122*M231+ M221*M132-M222*M131},
      {-(M112*M233-M113*M232+M212*M133-M213*M132),
       M111*M233-M113*M231+   M211*M133-M213*M131,
       -(M111*M232-M112*M231+ M211*M132-M212*M131)},
      { M112*M223-M113*M222+ M212*M123-M213*M122,
       -(M111*M223-M113*M221+ M211*M123-M213*M121),
         M111*M222-M112*M221+M211*M122-M212*M121}}) 


-- wedge2:     input:   single 3x2 matrices u, 
--             output:  single 3x1 matrix  \wedge^2 (u) 
wedge2 = (u) ->  
   matrix {{det(submatrix (u,{1,2},{0,1}))},
           {-det(submatrix (u,{0,2},{0,1}))},
           {det(submatrix (u,{0,1},{0,1}))}}

-- mixedwedge:   input two 3x2 matrices u1 and u2, output 
--     single 3x1 matrix      u1 \wedge u2 + u2 \wedge u1 
mixedwedge = (u1,u2) -> 
 matrix{
    {det (submatrix(u1,{1},{0,1}) || submatrix(u2,{2},{0,1}))+
     det (submatrix(u2,{1},{0,1}) || submatrix(u1,{2},{0,1}))},
    {- (det (submatrix(u1,{0},{0,1}) ||
             submatrix(u2,{2},{0,1}))+
        det (submatrix(u2,{0},{0,1}) ||
             submatrix(u1,{2},{0,1})))},
       {det (submatrix(u1,{0},{0,1}) ||
             submatrix(u2,{1},{0,1}))+
         det (submatrix(u2,{0},{0,1}) ||
              submatrix(u1,{1},{0,1}))}}
x =  matrix {{x1,x2,x3}}
y =  matrix {{y1,y2,y3}}
A1 = matrix {{0,A1p12,A1p13},{-A1p12,0,A1p23},{-A1p13,-A1p23,0}}
A2 =  matrix {{0,A2p12,A2p13},{-A2p12,0,A2p23},{-A2p13,-A2p23,0}}
B1 =  matrix {{B1p11,B1p12,B1p13},{B1p12,B1p22,B1p23},
      {B1p13,B1p23,B1p33}}
B2 =  matrix {{B2p11,B2p12,B2p13},{B2p12,B2p22,B2p23},
      {B2p13,B2p23,B2p33}}

-- A1,A2 are generic 3x3 skew,  B1,B2 are generic 3x3 symmetric

ID =  minors (2, matrix {{ y1,y2,y3, z*x1,z*x2, z*x3}, 
                         { x1,x2,x3, y1,y2, y3}})
  adjB1 = adjointMatrix (B1);
  adjB1p11 = adjB1_0_0; adjB1p22 = adjB1_1_1; adjB1p33 = adjB1_2_2;
  adjB1p12 = adjB1_1_0; adjB1p13 = adjB1_2_0; adjB1p23 = adjB1_2_1;
  adjB2 = adjointMatrix (B2);
  adjB2p11 = adjB2_0_0; adjB2p22 = adjB2_1_1; adjB2p33 = adjB2_2_2;
  adjB2p12 = adjB2_1_0; adjB2p13 = adjB2_2_0; adjB2p23 = adjB2_2_1;

a  =  det matrix {{x1,x2,x3},{A1p23,-A1p13,A1p12},
                    {A2p23,-A2p13,A2p12}};
invof2 = substitute ((1_S/2_S), S);

b = (invof2 * (y* B1* adjB2*B1*transpose(y)+
              y* B2* adjB1*B2*transpose(y)+
          2*x*A1* adjB2*(transpose A1)*transpose(x)+
          2*x*A2* adjB1*(transpose A2)*transpose(x)-
          2*x*A1*doubleAdjointMatrix (B1,B2)*
                 transpose(A2)*transpose(x)-
          z*x*B1*adjB2*B1*transpose(x)-
          z*x*B2*adjB1*B2*transpose(x)+
    4*y*B1*adjB2*transpose(A1)*transpose(x)+ 
    4*y*B2*adjB1*transpose(A2)*transpose(x)-
    (A1p12*(x1*y2-x2*y1)+A1p13*(x1*y3-x3*y1)+ 
                 A1p23*(x2*y3-x3*y2))*
    (B1p11* adjB2p11+B1p22* adjB2p22+ B1p33*adjB2p33 + 
     2*(B1p12* adjB2p12+B1p13* adjB2p13+ B1p23*adjB2p23))  -
    (A2p12*(x1*y2-x2*y1)+A2p13*(x1*y3-x3*y1)+ 
                 A2p23*(x2*y3-x3*y2))*
    (B2p11* adjB1p11+B2p22* adjB1p22+ B2p33*adjB1p33 +
    2*(B2p12* adjB1p12+B2p13* adjB1p13+B2p23*adjB1p23))))_0_0;

u1 = (-A1*(transpose x) + B1 *(transpose y)) | 
     (-A2*(transpose x) + B2 *(transpose y))
u2 = (-B1*(transpose x)) | (-B2*(transpose x))
sigma1 = wedge2 (u1) + z * wedge2 (u2)
sigma2 = mixedwedge (u1,u2)
f1 = (x*A1*transpose(y)+y*B1*transpose(y)-z*x*B1*transpose(x))_0_0;
f2 = (x*A2*transpose(y)+y*B2*transpose(y)-z*x*B2*transpose(x))_0_0;
l1 = z*x1*s0+y1*s1+sigma1_0_0;
l2 = z*x2*s0+y2*s1+sigma1_0_1;
l3 = z*x3*s0+y3*s1+sigma1_0_2;
l4 = y1*s0 + x1*s1+ sigma2_0_0;
l5 = y2*s0 + x2*s1+ sigma2_0_1;
l6 = y3*s0 + x3*s1+ sigma2_0_2;
q = s1^2-z*s0^2-a*s0+b; 

IY =ideal (f1,f2,l1,l2,l3,l4, l5,l6,q);

--  The following is an example of substitution to a nongeneric case
subslist = { x1 => x1,x2 => 0,x3 => x2,
y1 => y1,y2 => y2,y3 => 0,  z => 5,
A1p12 => 1,A1p13=>2 ,A1p23=>0,
A2p12 => 1,A2p13=> 3 ,A2p23=> 1, 
B1p11 =>1, B1p12 => 5, B1p13=> 0, 
    B1p22 =>1, B1p23=>1 ,B1p33=> 1, 
B2p11=> -1,B2p12 => 2,B2p13=> 1,
   B2p22=>1,B2p23=>7 ,B2p33=>1}; 
specificIYbefore = substitute ( IY,subslist);
specificS = kk [  y1,y2,y3,x1,x2,x3,s0,s1, Degrees => {8:1}] 
specificIY = sub (specificIYbefore,specificS)
isHomogeneous specificIY  -- true 
codim specificIY
betti res specificIY             
\end{verbatim}

\bigskip
\noindent
Stavros Papadakis \\
Center for Mathematical Analysis, Geometry, and Dynamical Systems \\
Departamento de Matem\'atica,  Instituto Superior T\'ecnico \\
Universidade T\'ecnica de Lisboa \\
1049-001 Lisboa, Portugal \\
e-mail: papadak@math.ist.utl.pt

\end{document}